\documentclass[12pt]{article}
\usepackage{amsfonts,amsmath,amsthm,mathrsfs,amssymb,graphicx}
\usepackage{ dsfont}
\usepackage{stmaryrd}
\usepackage{appendix}
\usepackage{subfigure}
\usepackage{url}
\usepackage{hyperref}
\usepackage{multicol}
\usepackage{verbatim}
\usepackage{tikz}
\usepackage{hyperref}
\usepackage{enumerate}
\usepackage{wrapfig}

\usepackage[]{algorithm2e}

\newtheorem{theorem}{Theorem}[section]
\newtheorem{corollary}{Corollary}[section]
\newtheorem{definition}{Definition}[section]

\newtheorem{lemma}{Lemma}[section]

\def\proof{\noindent{\bf Proof: }}
\def\eproof{\hfill $\Box$ \par}

\makeatletter
\newcommand*{\rom}[1]{\expandafter\@slowromancap\romannumeral #1@}
\makeatother

\title
{Knot Morphing Algorithm for Quantum `Fragile Topology'}

\author{K. E. Jordan\thanks{IBM.
      Email:  {\tt kjordan@us.ibm.com.}} and
	J.Li\thanks{Department of Mathematics,
       University of Connecticut, 196 Auditorium Road, Unit 3009, Storrs, CT 06269.  
       Email:  {\tt ji.li@uconn.edu.}}
       and T. J. Peters\thanks{Department of Computer Science \& Engineering,
       University of Connecticut, Storrs, CT 06269.        
       Email:  {\tt tpeters@cse.uconn.edu.} This author acknowledges, with 
      appreciation, generous financial support from IBM Research, under award
       IBM-TJP-6328340.  All statements made, inclusive of any errors, by this
      author are the solely the responsibility of the author, not of IBM Research.} 
      }  
\date{\today}

\begin{document}
\vspace{-4ex}  
\maketitle

\vspace{-7ex}
\begin{abstract}
A knot theoretic algorithm is proposed to model `fragile topology' of quantum physics.
\end{abstract}
\vspace{-6ex}

\section{Introduction}
\label{sec:intro}

\begin{figure}[h!]
\centering
    \subfigure[Unknot]
   {
   \scalebox{0.8}{\includegraphics[height=5cm]{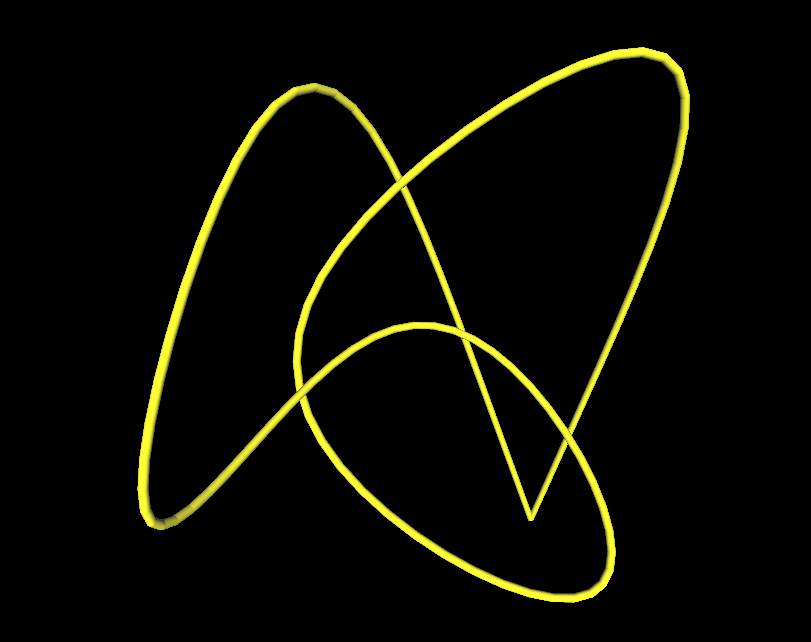}}
   \label{fig:end1}
   }
    \subfigure[Knot $4_1$]
   {
   \scalebox{0.8}{\includegraphics[height=5cm]{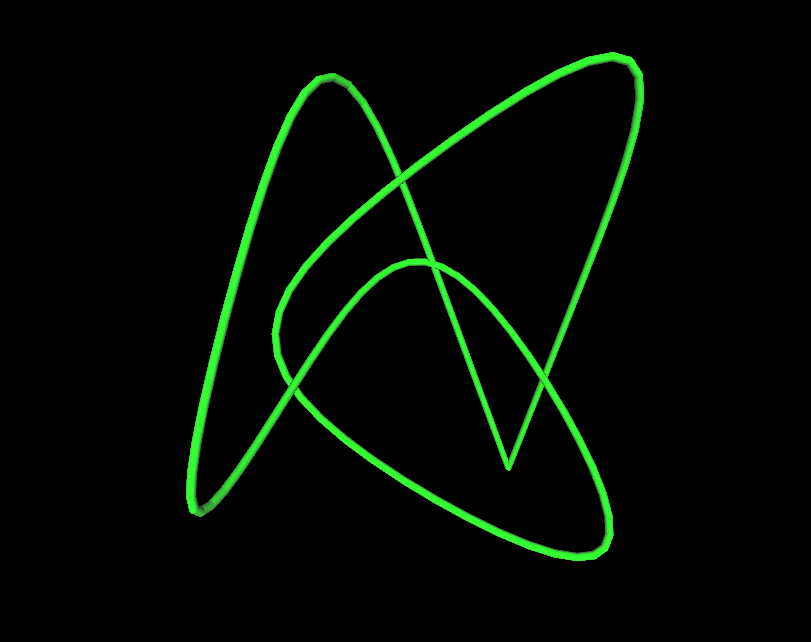}}
   \label{fig:end2}
   }
   \subfigure[Superimposed Knots]
   {
   \scalebox{0.8}{\includegraphics[height=5cm]{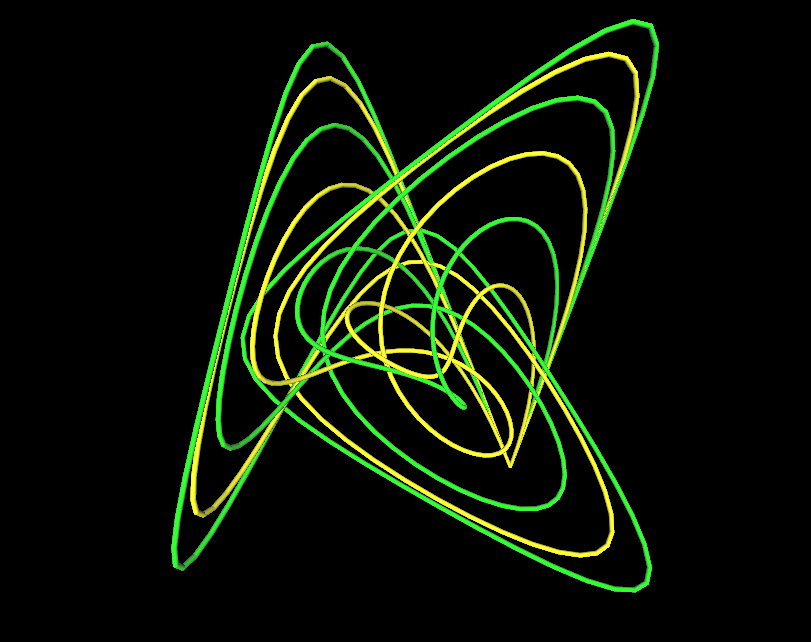}}
   \label{fig:ends6}
   }
\vspace{-1ex}
    \caption{Knots as Boundary Curves for Ruled Surfaces
}
    \label{fig:surfcon}
\end{figure}

Temporary, local topological changes characterize `fragile topology' in quantum physics.  As a mathematical model, envision the unkot of Figure~\ref{fig:end1} linearly sweeping out a nonselfintersecting surface, normal to the image. 
If one boundary curve is then morphed into $4_1$ of Figure~\ref{fig:end2} a selfintersection  appears in the surface.  The instantaneous change into a non-manifold mimics the spatiotemporal behavior of fragile topology. The overlaid images of Figure~\ref{fig:ends6} indicate the topological complexity possible by morphing multiple knots.   An algorithm enables computational experiments to rigorously analyze the topological instabilities of the quantum physics.

\pagebreak

Knots of one type can instantaneously deform into another type.  When modeled as B\'ezier curves, these knots can be embedded into surfaces and morphed to model ``fragile topology'', as appears in twisted bilayer graphine. 
An algorithm is presented to generate surfaces for analysis of these quantum phenomena.  That analysis may be complemented by emerging machine learning methods for high-throughput screening of material properties~\cite{sendek2017holistic}.

The presentation of related work ensues, followed by the underlying mathematics, the generative algorithm, then ending with conclusions and discussion of future work.

\section{Related Work}
\label{sec:relwork}

Loops are prominent curves in the theory of quantum materials~\cite{bradlyn2019disconnected,po2018fragile} and a loop is just the unknot.    Quantum material bands are also curves, with some so topologically unstable to be described as
``fragile topology''~\cite{po2018fragile}.  These curve models provide insights into twisted bilayer graphene as a fundamental quantum material~\cite{lian2018landau,po2018fragile,wieder2018axion}.  
In the work presented here, the unknot is morphed into a more complex knot to model transitions similar to fragile topology.

Knot theory has prompted discovery in quantum physics, as in three representative examples.
The concept of knotted electromagnetic fields led to the first experimental images of topological three-dimensional skyrmions~\cite{lee2018synthetic}.  The timing of knot creation~\cite{ollikainen2017quantum} and decay~\cite{ollikainen2019decay} are central to modeling quantum properties, as further explored here.  
Recent analysis of quantum manifolds~\cite{huang2017complexity} shares fundamental ideas with the surfaces presented here.
Historically, the literature on the role of knots for modeling energy and other physics phenomena spans decades~\cite{calvo2002physical}.

Many analyses of knots, particularly physical properties, proceed from creation of a tubular neighborhood~\cite{Hirsch} about the knot, which relies on a computation of its radius~\cite{cantarella2002minimum,cantarella2012shapes,DenneSullivan2008}.  A similar, but more restricted length is the basis for the ruled surfaces created here.

Standard, basic terminology for knot theory appears in the monorgraphs~\cite{Armstrong1983,Livingston1993,rolfsen2003knots}, with the first  also being a source for standard topological definitions.  Background for B\'ezier curves and ruled surfaces appears in the monograph~\cite{G.Farin1990}, within the broader class of splines and NURBS.
The paper~\cite{whited2010ball} defines different types of morphing on curves and analyzes their behaviors.

A previous result~\cite{marinelli2019exact} morphed an unknot into a trefoil, when both were represented as B\'ezier curves.  This was done by extensive computational experiments, essentially a Monte Carlo method, to discover how far to move a single control point.  The approach was similar to 
D-NURBS~\cite{qin1996d}, which were conceived to change spline geometry by integration with physics-based input.  The challenge over both approaches is to determine the perturbations of the control points.  These perturbations range over all vectors in $\mathbb{R}^3$, while the methods presented here restrict new control points to a given polygon.   Equivalence of a B\'ezier curve to a given piecewise linear knot has been established under a convergent sequence~\cite{li2015topological}, with sufficient~\cite{LPJZ16} and necessary conditions~\cite{jordan2019subdivision} for some knots, by contributions of two of the present authors.

\section{Mathematical Modeling by Ruled Surfaces}
\label{sec:mathmod}

A 2-manifold is created by sweeping the unknot of Figure~\ref{fig:end1} normal to the plane shown.  Although Figure~\ref{fig:end1} shows four crossings, these can be eliminated by Reidemeister moves~\cite{Livingston1993}. To be nonselfintersecting, the length of the sweep is chosen to be less than 1/2 of the minimum of the four distances from the pairs of preimages of the projected intersections.  A sweep is a special case of a ruled surface defined, below.

While these four crossings of Figure~\ref{fig:end1} are, visually, very similar to those depicted for $4_1$ of Figure~\ref{fig:end2}, the topological differences are significant, as these knots are not equivalent.    Attempting to modify this swept surface by locally morphing one boundary curve into $4_1$ will produce a selfintersecting surface.  While this phenomenon is readily understood within knot theory, the additional advantage provided here is algorithmic morphing to permit computational experiments for

 \begin{itemize}
\vspace{-0.09in}
\item localizing topological changes, as opposed to just proving their existence, and
\vspace{-0.09in}
\item mathematically modeling changes in quantum physics.
\vspace{-0.09in}
\end{itemize}

\subsection{Knots as Surface Boundaries}
\label{ssec:kboundaries}

The following definitions are standard, but included for the sake of completeness.
The equivalence relation from ambient isotopy is fundamental in knot theory.

\begin{definition}
Let $X$ and $Y$ be two subspaces of $\mathbb{R}^3$.  A continuous function
\[ H:\mathbb{R}^3 \times [0,1] \to \mathbb{R}^3 \]
is an {\bf ambient isotopy} between $X$ and $Y$ if $H$ satisfies the following conditions:
\begin{enumerate}
\vspace{-0.1in}
\item $H(\cdot, 0)$ is the identity,
\vspace{-0.1in}
\item $H(X,1) = Y$, and
\vspace{-0.1in}
\item $\forall t \in [0,1], H(\cdot,t)$ is a homeomorphism from
$\mathbb{R}^3$ onto $\mathbb{R}^3$.
\vspace{-0.1in}
\end{enumerate}
The sets $X$ and $Y$ are then said to be {\bf ambient isotopic}.
\label{def:aiso}
\end{definition}

All the knots considered for the scope of this work are assumed to be `tame'.  

\begin{definition}
\label{def:tame}
A {\bf knot is tame} if it is ambient isotopic to a knot  formed from finitely many linear segments.  
\end{definition}

For a tame knot, there always exists a planar
projection such that the double points occur exactly at the undercrossings and overcrossings.  Figure~\ref{fig:end2} is such a  projection of $4_1$.  While Figure~\ref{fig:end1} is \emph{not} such a projection of the unknot, it is shown for its visual similarity to $4_1$.

\begin{definition}
Let two curves be denoted by
\[c_1: [0,1] \rightarrow \mathbb{R}^3, \hspace{1ex} \mathrm{and} \hspace{1ex}  
c_2: [0,1] \rightarrow \mathbb{R}^3,\]
then, for $u \in [0,1]$ and $v \in [0,1]$, a {\bf ruled surface}, denoted by  $\Psi$, is defined by
\[ \Psi(u,v) = (1 - v)  c_1(u) + v  c_2 (u).\]
\label{def:ruled}
\end{definition}

\begin{lemma}
Let $c_1, c_2$ denote two tame knots, defining a ruled surface $S$.  If $S$ is nonselfintersecting, then $c_1$ and $c_2$ are ambient isotopic.  
\label{lem:ruled-iso}
\end{lemma}

\proof
Following Definition~\ref{def:aiso}, for any $s \in [0,1]$ and $\forall t \in [0,1],$ and for the given $c_1$ and $c_2$, let $H^*$ denote the function  
\[ H^*: [0,1]^2 \to \mathbb{R}^3, \]
\[H^*(s,t) = (1-t) c_1(s) + t c_2(s).\]  
Clearly, $H^*(s,t)$ is continuous.  Also, for each $t \in [0,1], H^*( \cdot, t)$ is now shown to be injective.  Suppose there exist $s_1, s_2 \in [0,1],$ with $s_1 \neq s_2$ such that $H^*(s_1,t) = H^*(s_2,t)$ -- a contradiction to $S$ being nonselfintersecting.  The function $H^*$ has a compact domain and can be extended to an ambient isotopy of compact support~\cite{andersson2000equivalence} in $\mathbb{R}^3$, so that $c_1$ and $c_2$ are ambient isotopic.
\eproof

\begin{corollary}
Let $c_1, c_2$ denote two tame knots of different knot type.  Then there does not exist any nonseflintersecting ruled surface between $c_1$ and $c_2$.
\label{cor:non-iso}
\end{corollary}

Despite the visual similarities between the unknot of Figure~\ref{fig:end1} and $4_1$ of Figure~\ref{fig:end2}, any ruled surface between them must necessarily have a selfintersection.  A more detailed analysis of the appearance of this selfintersection over time becomes possible by use of the curve morphing described in the next subsection.

\subsection{Morphing of B\'ezier Curves}
\label{ssec:morphbez}

\begin{definition}\label{def:c}
Let $\alpha(t)$ denote the {\bf parameterized B\'ezier curve} of degree $n$ with control points $P_i \in \mathbb{R}^3$, defined by
\vspace{-1ex}
$$
\alpha(t)=\sum_{i=0}^{n}{B_{i,n}(t)P_i}, \hspace{2ex}
t\in[0,1]
$$
where $B_{i,n}(t) = \left(\!\!\!
  \begin{array}{c}
	n \\
	i
  \end{array}
  \!\!\!\right)t^i(1-t)^{n-i}$.

The curve $\mathcal{P}$ formed by piecewise linear (PL) interpolation on the ordered list of points $P_0,P_1,\ldots,P_n$ is called the {\bf control polygon} and is a PL approximation of $\alpha$.
\end{definition}

Following previous work~\cite{li2015topological} on an initial control polygon $\mathcal{P}$, a sequence of new control polygons and B\'ezier curves will be generated.  Let ${\mathcal{P}}^{(0)} = \mathcal{P}$.  For $j \geq 0$, generate 
${\mathcal{P}}^{(j+1)}$ by the insertion of midpoints, as new control points, between all of the control points of ${\mathcal{P}}^{(j)}$ 

To avoid trivial and degenerate cases, for any control polygon $\mathcal{P}$, \emph{initially given as input}, it is assumed that for $i > 0$,
\begin{itemize}
\vspace{-0.1in}
\item $P_i \neq P_i*$, for any $i \neq i*$,
\vspace{-0.1in}
\item no more than two consecutive control points are collinear,
\vspace{-0.1in}
\item $n \geq 4$.
\vspace{-0.09in}
\end{itemize}

\vspace{1ex}

\begin{wrapfigure}{r}{0.5\textwidth}
  \vspace{-10ex}
  \begin{center}
    \includegraphics[width=0.40\textwidth]{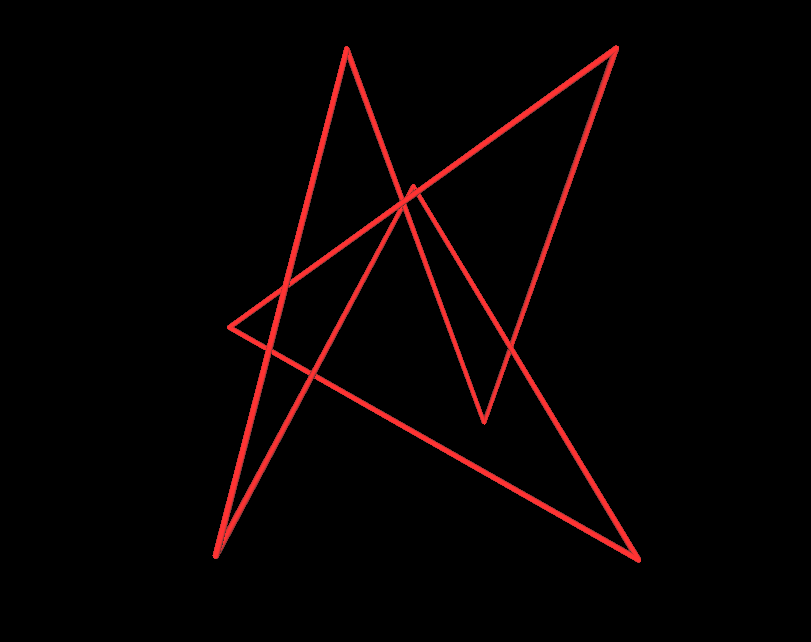}
  \end{center}
  \vspace{-4ex}
  \caption{0th Control Polygon}
  \label{fig:4_1-0th-poly}
    \vspace{-5ex}
\end{wrapfigure}
A control polygon, which is also the knot $4_1$, is shown in Figure~\ref{fig:4_1-0th-poly}

\begin{theorem}
For iterative insertion of midpoints, the distance between the control polygon and the generated B\'ezier curves converges to zero.
\label{thm:colli-conv}
\end{theorem}
Moreover, after sufficiently many iterations, the B\'ezier curves are ambient isotopic to the initial control polygon~\cite{jordan2019subdivision,LPJZ16,li2015topological}.  It has been shown~\cite{jordan2019subdivision} that iterations 0 through 3 produce unknotted B\'ezier curves, while the 4th iteration results in the B\'ezier curve being $4_1$, as shown in Figure~\ref{fig:end2}.  (All the control polygons formed are ambient isotopic, trivially, by the identity.)

\section{Computational and Algorithmic Framework}
\label{sec:val}

The topological change over a ruled surface has been depicted over two static images of  Figures~\ref{fig:end1} and \ref{fig:end2}.   The first ruled surface is a manifold and the second is a non-manifold.  The key question regarding fragile topology is
\begin{center}
Can the moment in time be modeled when this change occurs?
\end{center}

The underlying mathematics afford opportunities for effective visual and numerical methods to approximate the timing of topological change.  

For more complex models of fragile topology, consider Figure~\ref{fig:ends6}, spanning 6 knots, 4 of them being the unknot and two being $4_1$.  These correspond to the 0th to 5th iteration of midpoint insertions.  These are similar to frames in a movie.  The animation community has developed many `in-betweening' techniques to produce frames between isolated images~\cite{di2001automatic}, as tools to be considered for the proposed temporal analysis.

Skinned surfaces are generalizations of ruled surfaces, extending beyond two boundary curves to interpolating over arbitrarily many curves. By algorithmically morphing the interpolated curves, the generation of selfintersections over time could be investigated.  
Also, a surface could be incrementally built over these 6 knots, by starting at the first two and then successively adding ruled surfaces.  Curve morphing would then support spatiotemporal analysis of selfintersections.

While the insertion of new control points has been presented at all midpoints, simultaneously, those two constraints can both be relaxed. Specifically,
\begin{itemize}
\vspace{-0.09in}
\item Any number of control points, even just one, can be added at any iteration.
\vspace{-0.09in}
\item The added points need not be midpoints, just collinear with existing control points.
\vspace{-0.09in}
\end{itemize}
Exercising the first generalization would provide finer control over time.  Use of the second would produce more local geometric control.  For example, if a local region of topological interest is highly effected by insertion of a particular control point, then more control points can be inserted nearby to gain further understanding.  Publicly available data bases of piecewise linear knots (also known as stick knots)~\cite{KServer} provide a rich store of examples for extensive computational experiments for modeling fragile topology.  Algorithm~\ref{algo:fragile} is central to these proposed experiments.

\begin{algorithm}[h]
\KwData{K = PL knot as a control polygon.   \hspace{9ex} $\sharp$ Knot type of K is known.}  
 \KwResult{A ruled surface to model fragile topology.}  
initialization to K \hspace{33ex}  $\sharp$ Sweep a nonselfintersecting surface.\\
\While{Morphing one B\'ezier boundary curve} 
{Visualize emerging ruled surfaces   \hspace{14ex} $\sharp$  Experiment with adding control points.\\
\eIf{Selfintersection is determined}{
 Dynamically visualize ruled surface. \hspace{8ex} $\sharp$  Experiment with fragile topology changes. \;
}{
Create alternative morphs.  \hspace{17ex} $\sharp$  Modify control points on one boundary. \;\;
}
}
\caption{Modeling Fragile Topology with Surfaces}
\label{algo:fragile}
\end{algorithm}

The advantage of these topology modifications is contrasted with a previous example~\cite{marinelli2019exact} of an unknotted control polygon with an unknotted B\'ezier curve.  Extensive experiments were done to modify the control points to change the knot type of the modified B\'ezier curve.  While this was ultimately successful, it was highly empirical, providing little future guidance to 
narrow the infinitely many options available.  Similar complexity of shape dependencies has been echoed elsewhere~\cite{andersson2000equivalence,qin1996d,sequinbeauty}.

The constraints imposed here restrict the possibilities of new control points so that more effective assistive algorithms are possible.  This can permit massive generation of surface examples for machine learning evaluations from instantiations of piecewise linear knots~\cite{rawdon2002upper} of optimal surface properties for preferred physics from fragile topology.  

An initial software prototype was created in Python and was used to generate the images here.  This was done at a `proof of concept' level to refine the mathematics presented.   However, the present user interface should be enhanced
for broad experiments of many data sets or for interactive changes.  The existing code points to broader computational experiments for quantum phsyics. 

\section{Conclusions and Future Work}
\label{sec:concfu}

Knots are boundaries for surfaces constructed to model the occurrence of `fragile topology' of quantum properties of materials, such as twisted bilayer graphene.  The key topological property captured is an instantaneous appearance of a local topological change.  This appears to be the first application of a sequence of converging B\'ezier curves to quantum physics, inclusive of an algorithm for further computational experiments to prompt insight.

Four specific questions are posed for specific future work and software design and implementation, relative to the morphing process presented:
\begin{itemize}
\vspace{-0.09in}
\item Can parameters be established to let the physics perturb candidate surfaces?
\vspace{-0.09in}
\item  Does further modeling value arise from extending to skinned surfaces?
\vspace{-0.09in}
\item Are there manageable hypotheses to constrain selfintersections to interior of surfaces?
\vspace{-0.09in}
\item What is an optimal mix of numerical methods and visual analytics for intersection detection?
\vspace{-0.09in}
\end{itemize}
\vspace{-1ex}

The opportunity is to integrate knot theory, spline geometry and quantum physics to algorithmically create multiple candidate surfaces, inclusive of ones not yet even imagined -- permitting computational experiments to discover nuances of fragile topology within quantum physics.

\section{Acknowledgments}
\label{sec:ack}

The authors, acknowledge, with appreciation, informative discussions on materials science and machine learning with Q. Yang of the University of Connecticut, as well as her insightful and constructive comments on a draft of this paper.

\small{
\bibliographystyle{plain}
\bibliography{TJP-vis-knot.bib}

\begin{thebibliography}{10}

\bibitem{andersson2000equivalence}
L-E Andersson, T.~J. Peters, and N.~F. Stewart.
\newblock Equivalence of topological form for curvilinear geometric objects.
\newblock {\em International Journal of Computational Geometry \&
  Applications}, 10(06):609--622, 2000.

\bibitem{Armstrong1983}
M.~A. Armstrong.
\newblock {\em Basic Topology}.
\newblock Springer, New York, 1983.

\bibitem{bradlyn2019disconnected}
B.~Bradlyn, Z.~Wang, J.~Cano, and B.~A. Bernevig.
\newblock Disconnected elementary band representations, fragile topology, and
  {W}ilson loops as topological indices: An example on the triangular lattice.
\newblock {\em Physical Review B}, 99(4):045140, 2019.

\bibitem{calvo2002physical}
J.~A. Calvo, K.~C. Millett, and E.~J. Rawdon.
\newblock {\em Physical Knots: Knotting, Linking, and Folding Geometric Objects
  in $\mathbb{R}^3$}, volume 304.
\newblock American Mathematical Society, 2002.

\bibitem{cantarella2002minimum}
J.~Cantarella, R.~B Kusner, and J.~M. Sullivan.
\newblock On the minimum ropelength of knots and links.
\newblock {\em Inventiones Mathematicae}, 150(2):257--286, 2002.

\bibitem{cantarella2012shapes}
J.~Cantarella, A.~LaPointe, and E.~J. Rawdon.
\newblock Shapes of tight composite knots.
\newblock {\em Journal of Physics A: Mathematical and Theoretical},
  45(22):225202, 2012.

\bibitem{DenneSullivan2008}
E.~Denne and J.~M. Sullivan.
\newblock Convergence and isotopy type for graphs of finite total curvature.
\newblock In A.~I. Bobenko, J.~M. Sullivan, P.~Schr{\"o}der, and G.~M. Ziegler,
  editors, {\em Discrete Differential Geometry}, pages 163--174. Birkh{\"a}user
  Basel, 2008.

\bibitem{di2001automatic}
F.~Di~Fiore, P.~Schaeken, K.~Elens, and F.~Van~Reeth.
\newblock Automatic in-betweening in computer assisted animation by exploiting
  2.5 d modelling techniques.
\newblock In {\em Proceedings Computer Animation 2001. Fourteenth Conference on
  Computer Animation (Cat. No. 01TH8596)}, pages 192--200. IEEE, 2001.

\bibitem{G.Farin1990}
G.~E. Farin.
\newblock {\em Curves and Surfaces for Computer-Aided Geometric Design: A
  Practical Guide}.
\newblock Academic Press, Inc., 1996.

\bibitem{Hirsch}
M.~W. Hirsch.
\newblock {\em Differential Topology}.
\newblock Springer, New York, 1976.

\bibitem{huang2017complexity}
Z.~Huang and A.~V. Balatsky.
\newblock Complexity and geometry of quantum state manifolds.
\newblock {\em arXiv preprint arXiv:1711.10471}, 2017.

\bibitem{jordan2019subdivision}
K.~E. Jordan, K.~Marinelli, T.~J. Peters, J.~A. Roulier, and P.~Zaffetti.
\newblock Subdivision of {B}{\'e}zier curves for ambient isotopy in molecular
  modeling.
\newblock {\em Topology and its Applications}, 259:311--322, 2019.

\bibitem{lee2018synthetic}
W.~Lee, A.~H. Gheorghe, K.~Tiurev, T.~Ollikainen, M.~M{\"o}tt{\"o}nen, and
  D.~S. Hall.
\newblock Synthetic electromagnetic knot in a three-dimensional skyrmion.
\newblock {\em Science Advances}, 4(3):eaao3820, 2018.

\bibitem{LPJZ16}
J.~Li, T.~J. Peters, K.~E. Jordan, and P.~Zaffetti.
\newblock Computational topology: Isotopic convergence to a stick knot.
\newblock {\em Topology and its Applications}, 206:276--283, 2016.

\bibitem{li2015topological}
J.~Li, T.~J. Peters, K.~Marinelli, E.~Kovalev, and K.~E. Jordan.
\newblock Topological subtleties for molecular movies.
\newblock {\em Topology and its Applications}, 188:91--96, 2015.

\bibitem{lian2018landau}
B.~Lian, F.~Xie, and B.~A. Bernevig.
\newblock The {L}andau level of fragile topology.
\newblock {\em arXiv preprint arXiv:1811.11786}, 2018.

\bibitem{Livingston1993}
C.~Livingston.
\newblock {\em Knot Theory}, volume~24 of {\em Carus Mathematical Monographs}.
\newblock Mathematical Association of America, Washington, DC, 1993.

\bibitem{marinelli2019exact}
K.~Marinelli and T.~J. Peters.
\newblock Exact computation for existence of a knot counterexample.
\newblock {\em Applied General Topology}, 20(1):251--264, 2019.

\bibitem{ollikainen2019decay}
T.~Ollikainen, A.~Blinova, M.~M{\"o}tt{\"o}nen, and D.~S. Hall.
\newblock Decay of a quantum knot.
\newblock {\em Physical Review Letters}, 123(16):163003, 2019.

\bibitem{ollikainen2017quantum}
T.~Ollikainen, S.~Masuda, M.~M{\"o}tt{\"o}nen, and M.~Nakahara.
\newblock Quantum knots in {B}ose-{E}instein condensates created by
  counterdiabatic control.
\newblock {\em Physical Review A}, 96(6):063609, 2017.

\bibitem{po2018fragile}
H.~C. Po, H.~Watanabe, and A.~Vishwanath.
\newblock Fragile topology and {W}annier obstructions.
\newblock {\em Physical review letters}, 121(12):126402, 2018.

\bibitem{qin1996d}
H.~Qin and D.~Terzopoulos.
\newblock D-{NURBS}: a physics-based framework for geometric design.
\newblock {\em IEEE Transactions on Visualization and Computer Graphics},
  2(1):85--96, 1996.

\bibitem{rawdon2002upper}
E.~J. Rawdon and R.~G. Scharein.
\newblock Upper bounds for equilateral stick numbers.
\newblock {\em Contemporary Mathematics}, 304:55--76, 2002.

\bibitem{rolfsen2003knots}
D.~Rolfsen.
\newblock {\em Knots and Links}, volume 346.
\newblock American Mathematical Soc., 2003.

\bibitem{sendek2017holistic}
Austin~D Sendek, Qian Yang, Ekin~D Cubuk, Karel-Alexander~N Duerloo, Yi~Cui,
  and Evan~J Reed.
\newblock Holistic computational structure screening of more than 12000
  candidates for solid lithium-ion conductor materials.
\newblock {\em Energy \& Environmental Science}, 10(1):306--320, 2017.

\bibitem{sequinbeauty}
C.~H. Sequin.
\newblock The beauty of knots.
\newblock \\ \mbox{{\url
  {https://pdfs.semanticscholar.org/94e5/eb1b7d5f6a0afbdef36ab7cb5b5c194ed6cb.pdf}}}.

\bibitem{KServer}
Author Unknown.
\newblock The {K}not {S}erver.
\newblock \url{http://www.colab.sfu.ca/KnotPlot/KnotServer/ }, 2003.

\bibitem{whited2010ball}
B.~Whited and J.~Rossignac.
\newblock Ball-morph: Definition, implementation, and comparative evaluation.
\newblock {\em IEEE Transactions on Visualization and Computer Graphics},
  17(6):757--769, 2010.

\bibitem{wieder2018axion}
B.~J. Wieder and B.~A. Bernevig.
\newblock The axion insulator as a pump of fragile topology.
\newblock {\em arXiv preprint arXiv:1810.02373}, 2018.

\end{thebibliography}
}
\end{document}